\documentclass{elsarticle}
\usepackage{amsmath}
\usepackage{amsfonts}
\usepackage{amssymb}
\usepackage{graphicx}
\usepackage{xfrac}
\usepackage{color}
\usepackage{epsfig}
\usepackage{algorithm, algpseudocode}

\usepackage[utf8]{inputenc} 
\usepackage{tikz}
\usetikzlibrary{arrows,decorations.pathmorphing,backgrounds,positioning,fit,matrix, shapes, chains}
\tikzstyle{line} = [draw, -latex']
\usepackage[T1]{fontenc} 
\usepackage{mathtools}

\DeclarePairedDelimiter\floor{\lfloor}{\rfloor}

\newproof{pf}{Proof}

\newcommand{\Norm}[1]{\left\|#1\right\|}
\newcommand{\abs}[1]{\left|#1\right|}

\newcommand{\Null}{\operatorname{Null}}

\newcommand{\bigO}{\mathcal{O}}
\newcommand{\half}{\frac{1}{2}}
\newcommand{\R}{\mathbb{R}}

\newcommand{\Pol}{\mathbb{P}}
\newcommand{\xb}{\mathbf{x}}
\newcommand{\yb}{\mathbf{y}}
\newcommand{\zb}{\mathbf{z}}
\newcommand{\bb}{\mathbf{b}}
\newcommand{\gb}{\mathbf{g}}
\newcommand{\fb}{\mathbf{f}}
\newcommand{\mb}{\mathbf{m}}
\newcommand{\wb}{\mathbf{w}}

\newcommand{\lambdab}{\boldsymbol{\lambda}}
\newcommand{\etab}{\boldsymbol{\eta}}
\newcommand{\Phib}{\boldsymbol{\Phi}}

\newcommand{\iopt}{i_{\scriptsize \mbox{opt}}}
\newcommand{\nopt}{n_{\scriptsize \mbox{opt}}}
\newcommand{\ntensor}{n_{\scriptsize \mbox{tp}}}
\newcommand{\nelim}{n_{\scriptsize \mbox{elim}}}
\renewcommand{\tilde}{\widetilde}

\DeclareMathOperator*{\argmin}{argmin}
\DeclareMathOperator*{\argmax}{argmax}

\begin{document}
\begin{frontmatter}
\title{A Node Elimination Algorithm for Cubature of High-Dimensional Polytopes\tnoteref{t1}}
\tnotetext[t1]{This material is based upon work supported by the National
  Science Foundation under grant DMS-17207431.}

\author[smu]{Arkadijs Slobodkins} 
\ead{aslobodkins@smu.edu}

\author[smu]{Johannes Tausch\corref{cor1}}
\ead{tausch@smu.edu}
\cortext[cor1]{Corresponding author}

\address[smu]{Department of Mathematics, Southern Methodist University,
  Dallas, TX 75275, USA}

\begin{abstract}
  Node elimination is a numerical approach to obtain cubature rules
  for the approximation of multivariate integrals.  Beginning with a
  known cubature rule, nodes are selected for elimination, and a new,
  more efficient rule is constructed by iteratively solving the moment
  equations. This paper introduces a new criterion for selecting which
  nodes to eliminate that is based on a linearization of the moment
  equation. In addition, a penalized iterative solver is introduced,
  that ensures that weights are positive and nodes are inside the
  integration domain. A strategy for constructing an initial quadrature
  rule for various polytopes in several space dimensions is
  described. High efficiency rules are presented for two, three and
  four dimensional polytopes.
  The new rules are compared with rules that are obtained
  by combining tensor products of one dimensional quadrature rules and
  domain transformations, as well as with known analytically
  constructed cubature rules.
\end{abstract}

\end{frontmatter}

\section{Introduction}
Let $\Omega$ be a domain in $\R^d$. The goal of the construction of
cubature rules is to determine the
nodes $x_k$ and weights  $w_k$ in the cubature rule
\begin{equation}\label{cubaturerule}
  \int_{\Omega} \phi(x) w(x) \, dx \approx
  \sum_{k=1}^n \phi(x_k) w_k,
\end{equation}
such that the rule is exact for multivariate polynomials of
degree up to $p$
\begin{equation*}
  \Pol^d_p = \mbox{span} \left\{
    x^\alpha : \alpha_1 + \dots +\alpha_d \leq p \right\}.
\end{equation*}
Here, $\alpha$ is a multi index and  $x^\alpha =
x_1^{\alpha_1}\cdots x_d^{\alpha_d}$. The dimension of the linear
space $\Pol^d_p$ is 
\begin{equation*}
M=\mbox{dim }\Pol^d_p = \binom{p+d}{d}.
\end{equation*}
From the viewpoint of numerics the monomial basis is ill-conditioned,
and so we consider a more general basis $\phi_1,\dots, \phi_M$ of
$\Pol^d_p$, for instance, orthogonal polynomials. We set
\begin{equation}\label{def:vandermonde}
  \Phi(x) = \left[
    \begin{array}{c}\phi_1(x)\\ \vdots \\\phi_M(x) \end{array}
  \right], \qquad
  \bb = \left[
    \begin{array}{c}\int_\Omega \phi_1(x) w(x)\,dx\\ \vdots
      \\\int_\Omega \phi_M(x) w(x) \,dx \end{array}\right].
\end{equation}
Further, we write for the vector of nodes
and weights, $\xb\in\R^{dn}$, $\wb \in \R^n$, respectively, and
\begin{equation}\label{vandermonde}
\Phib(\xb) = \big[\Phi(x_1),\dots, \Phi(x_n) \big] \in \R^{M\times n}.
\end{equation}
Exactness in $\Pol^d_p$ means that the nodes and weights must be solutions
of moment equations
\begin{equation}\label{def:polysys}
\fb(\xb,\wb) = \Phib(\xb) \wb - \bb =
\mathbf{0}.
\end{equation}
which is a polynomial system in $N=(d+1)n$ unknowns and $M$ equations. 
We write \eqref{def:polysys} in more convenient form by combining nodes
and weights into one vector. Thus we let 
\begin{equation}\label{def:z}
z_k = [x_k, w_k] \in \R^{d+1}\qquad\mbox{and}\qquad  
\zb = [z_1, \dots, z_n]^T \in \R^N. 
\end{equation}
Since we are looking for nodes
in the domain $\Omega$ that have positive weights, a feasible cubature
rule is in the set
\begin{equation*}
  Z_n = \big\{ \zb\in\R^N :\, \fb(\zb) = \mathbf{0},
  \; x_k \in \Omega, \; w_k\geq 0, \; 1\leq k\leq n\big \}.
\end{equation*}
A solution of \eqref{def:polysys} in $Z_n$ is said to have quality PI,
which means that all weights are positive and all nodes are inside the
domain.  The construction of cubature rules has considerable
interest in several application areas. For instance, rules for three
dimensional polytopes are used in the $p$-version of the finite
element method~\cite{babuska88}. The Galerkin boundary element method
involves integrals over four dimensional
polytopes~\cite{sauter-schwab11,tausch22,manson-tausch19}. Likewise,
variational formulations of elliptic fractional PDEs in two or three
dimensions lead to four or six-dimensional
integrals~\cite{bonito-lei-pasciak19}.

Since \eqref{def:polysys} is a polynomial system, it can be approached
from the perspective of algebraic geometry, see ~\cite{cools97} and
the reference therein. In principle, the system can be solved using 
Gr\"obner bases and Buchberger's algorithm. However, the computational
complexity of these algorithms grows very quickly with the number of
variables, so that in practice the method can only provide solutions
for relatively low degree and dimension.  In addition, surprisingly
little is known about the solvability of \eqref{def:polysys}. In
particular, the smallest value of $n$ is for which $Z_n$ is non-empty
is generally unknown.

One can expect that the solutions of \eqref{def:polysys} live on
higher dimensional algebraic varieties if the system is
underdetermined, i.e., $N$ is greater than $M$.  If $N = M$ the
solutions are isolated points. This however, does not preclude that
there are rules with fewer nodes. We call a rule optimal, if $n$ is
the smallest integer such that $N \geq M$. For a collection of
cubature rules that were derived using analytic means we refer to
\cite{cools03}.

To illustrate the issues with multidimensional cubature, consider the case that
$\Omega$ is a hypercube. In this case cubature rules can be easily
found by forming tensor products of one-dimensional Gauss quadrature
rules. The resulting rules are solutions of \eqref{def:polysys}, but
these rules become highly underdetermined as the degree and the
dimension is increased.

For larger values of $p$ and $d$, analytical methods for constructing
solutions become very difficult, if not impossible, and therefore
methods that solve \eqref{def:polysys} with Newton-like methods and/or
stochastic methods have been researched extensively in the recent
years, see, e.g., \cite{jaskowiec-sukumar19,sudhakar-etal17,
  gentile-sommariva-vianello11, keshavarzzadeh-kirby-narayan18,
  vioreanu-rokhlin14}.

In this article we focus on generating quadrature rules by node
elimination, a technique that was initiated in the paper by Xiao and
Gimbutas~\cite{xiao-gimbutas10}. The idea is to begin with a known
overdetermined rule and select a node for elimination from the rule. The reduced
rule does not satisfy the moment equations, but is used as an initial
guess for solving this nonlinear system by the Gauss-Newton
iteration. The procedure is repeated until no further nodes can be
eliminated. 

Since nonlinear solvers typically converge only locally, the success
of such an approach strongly depends on how a node is selected for elimination. Xiao
and Gimbutas choose the node whose corresponding column in
\eqref{vandermonde} has the smallest Euclidean norm (scaled by the
weight).  In this paper we present a new elimination criterion which is based
on the linearization of the function $\fb(\xb,\wb)$ with the goal to
guarantee a close initial guess for the nonlinear solver. Thus our
node elimination procedure can be viewed as a predictor-corrector type
method that is used in path following algorithms, see~\cite{allg-geo03}.

Another important aspect of node elimination is that a solution is
only a useful cubature rule if the nodes are in $\Omega$ and the
weights are positive. In this article we describe a new penalized
optimization method that guarantees that nodes remain in the domain
and have positive weights. 

The outline of the remainder of this paper is as follows. In
section~\ref{sec:lsn} we describe the corrector step, which is a
constrained least squares Newton method. In section~\ref{sec:nodeelim}
we discuss the improved node elimination step (predictor). One
important aspect of node elimination is to come up with a good initial
cubature rule. In section~\ref{sec:initialguess} we describe how this can be
done for a variety of domains $\Omega$. Finally, in section
\ref{sec:numresult}, we present some cubature rules that we found
with the method. The codes of our node elimination
scheme and some of the cubature rules that we obtained can be found
on Github ~\cite{slobodkins}.

\section{Penalized Least Squares Newton Method}\label{sec:lsn}
As we already mentioned, the node elimination procedure is a
predictor-corrector type method. In this section we describe the
corrector step, which, for a point $\tilde \zb \not\in Z_n$, attempts
to find a nearby point on the solution manifold, i.e., $\bar \zb \in Z_n$. 

Since the goal is to obtain quadrature rules with nodes in the domain
that have positive weights, we enforce these constraints by adding a penalty
term. With $\zb$ defined as in \eqref{def:z}, we set
\begin{equation}\label{def:penalty:fcn}
\phi_\Omega(\zb) = \sum_{j=1}^n \left[ \varphi_\Omega( x_j ) + 
  \log\left( \frac{1}{w_j} \right) \right],
\end{equation}
where $\varphi_\Omega(\cdot)$ is a
function that is smooth in $\Omega$ and 
has a logarithmic singularity on the boundary of $\Omega$. For
instance, if $\Omega$ is a convex polytope given by the linear 
inequalities, then we set
\begin{equation*}
\Omega = \left\{ x\in \R^d: A x \leq b \right\} 
\quad\Rightarrow\quad \varphi_\Omega(x) = \sum_{\ell=1}^{\ell_A}
\log\left( \frac{1}{b_\ell - a_\ell^T x} \right),
\end{equation*}
where $a^T_\ell$ are the rows of $A$ and $\ell_A$ is the number of rows.
For the case of a sphere the penalty term is
\begin{equation*}
\Omega = \left\{ x\in \R^d: \Norm{x}_2 \leq 1 \right\} 
\quad \Rightarrow\quad \varphi_\Omega(x) = 
\log\left( \frac{1}{1-\Norm{x}_2} \right).
\end{equation*}

To derive an iterative solver that maps $\tilde \zb \not\in Z_n$ to
$\bar \zb \in Z_n$ let $\Delta \zb = \tilde \zb - \bar \zb$, and
consider the constrained optimization problem
\begin{equation}\label{nl:opt:problem}
\min\limits_{ \Delta \zb} \left\{
\half \Norm{\Delta \zb}^2 + t \,\phi_\Omega(\tilde \zb + \Delta \zb),\;
\tilde \zb + \Delta \zb \in Z_n \right\}.
\end{equation}
Here, the parameter $t\geq 0$
controls the strength of the penalty term. We will provide more
detail about how to determine $t$ later on.

Now linearize \eqref{nl:opt:problem} as follows
\begin{equation}\label{lin:opt:problem}
\begin{aligned}
\min \;       &\half \Norm{\Delta \zb}^2 + t\, \gb^T\Delta \zb \\
\Delta \zb:\; & \fb(\tilde \zb) + J \Delta \zb =  \mathbf{0},
\end{aligned}
\end{equation}
where $\fb:=\fb(\tilde \zb)\in \R^M$, $J := D\fb(\tilde \zb)\in \R^{M\times N}$
is the Jacobian of $\fb(\tilde \zb)$ and $\gb := \nabla\phi_\Omega(\tilde \zb) \in \R^N$ is the
gradient of the penalty term. Using Lagrange multipliers it follows that
\eqref{lin:opt:problem} is equivalent to the linear system
\begin{equation*}
\begin{aligned}
\Delta \zb + J^T \lambdab &= -t \gb \\
J \Delta \zb & = -\fb
\end{aligned}
\end{equation*}
where $\lambdab \in \R^M$ is the Lagrange multiplier. Eliminating
$\lambdab$ gives the solution of \eqref{lin:opt:problem} in the
the normal equations form
\begin{equation}\label{delta:z:ne}
\begin{aligned}
  \Delta \zb &= - J^T \left( J J^T\right)^{-1} \fb
           - t \left( I - J^T \left( J J^T\right)^{-1} J \right) \gb\\
             &= \Delta \zb^f + t  \Delta \zb^g.
\end{aligned}
\end{equation}
Here, $\Delta \zb^f = -J^T \left( J J^T\right)^{-1} \fb$ is the least squares solution of the
underdetermined system $\fb + J\Delta \zb = \mathbf{0}$, i.e., the
solution of \eqref{nl:opt:problem} when $t=0$. The second term 
$\Delta \zb^g = -\left( I - J^T \left( J J^T\right)^{-1} J \right) \gb$
is the orthogonal projection of $-\gb$ onto the nullspace of
$J$. The form in \eqref{delta:z:ne} makes clear that the solution of
\eqref{lin:opt:problem} for any parameter $t$ can be obtained by
computing $\Delta \zb^f$ and $\Delta \zb^g$ independently of $t$ and then forming the
appropriate linear combination.

For numerical purposes it is better to compute these two vectors with
the LQ factorization. If $J=LQ$ is the economy size factorization,
i.e., $L \in \R^{M\times M}$ is lower triangular and $Q \in
\R^{M\times N}$ has orthonormal rows, then
\begin{equation}\label{delta:zfg}
\begin{aligned}
  \Delta \zb^f &= -Q^T L^{-1} \fb,\\
  \Delta \zb^g &= -(I - Q^T Q) \gb.
\end{aligned}
\end{equation}
Once $\Delta \zb^f$ and $\Delta \zb^g$ have been computed, the
Newton update is
\begin{equation}\label{newton:update}
\tilde \zb  \leftarrow \tilde \zb + \Delta \zb^f + t \Delta \zb^g.
\end{equation}
We now describe how the 
parameter $t$ is determined. It is determined such that 
the nodes of next iterate have
maximal distance from the boundary and 
weights are as large as possible. We focus on the case that $\Omega$
is a convex polytope given by the linear inequalities $Ax\leq b$.

Substitution of the $j$-th node and weight of the Newton update
\eqref{newton:update} into the linear constraints results in
two types of inequalities
\begin{equation*}
  \begin{aligned}
b_\ell - a_\ell^T (\tilde x_j + \Delta x_j^f ) - t a_\ell^T \Delta
    x_j^g &\leq 0\\
-\tilde w_j +  \Delta w_j^f + t \Delta w_j^g &\leq 0
    \end{aligned}
\end{equation*}
for $1\leq j\leq n$, $1\leq\ell\leq \ell_A$. Here $\tilde x_j$,
$\Delta x_j^f$ and $\Delta x_j^g$ are the nodal components and
$\tilde w_j$, $\Delta w_j^f$ and $\Delta w_j^g$ are the weights of the
vectors $\tilde \zb$, $\Delta \zb^f$ and $\Delta \zb^g$,
respectively. We write these conditions collectively as
\begin{equation}\label{family}
m_k(t) = \beta_k + t \alpha_k \leq 0, \quad 1\leq k\leq (\ell_A+1)n,
\end{equation}
where the $\alpha_k$'s  list all inner products $a_\ell^T \Delta
x_j^g$ and all $\Delta w_j^f$'s and the $\beta_k$'s are defined
analogously.

Algorithm \ref{algo:minmax:family} finds the value of $t$ that
minimizes the maximum over $k$ for this 
family of straight lines. The algorithm starts with the maximum
at $t=0$ and follows the largest line until it intersects with another
line with positive slope.

\begin{algorithm}
\caption{Determines the min max of the family of straight
  lines in \eqref{family}. Here, $t_{jk}$ denotes the
intersection point of line $m_j$ with $m_k$.}
\label{algo:minmax:family}
\begin{algorithmic}
\State $t=0$
\State $k = \argmax\limits_{j} \beta_j$  
\While{$\alpha_k < 0$}
\State Find $k^* =  \argmin\limits_{k} \left\{ t_{jk},\; t_{jk}> t\right\}$
\State $t = t_{k^*k}$
\State $k^* = k$
\EndWhile
\end{algorithmic}
\end{algorithm}

If successful, the algorithm returns the $t$-value for which the min
max of the family of straight lines is achieved. This guarantees that
the next iterate has maximal distance from the boundary.  Note that
even though $\tilde \zb$ satisfies the constraints, the least squares
solution $\tilde \zb + \Delta \zb^f$ may not. Therefore, it is
possible that the returned value of $m_k(t)$ is positive in which case
at least one node does not satisfy the constraints.  In the latter case
one can resort to a damped Newton method, where $\Delta \zb^f$ is
multiplied by a small factor.  The theoretical cost of this algorithm
is $O(\ell_A^2 n^2)$. However, in practice, the while loop terminates
after a small number of steps, thus its computational cost is closer
to $O(\ell_A n)$ and negligible compared to the cost of computing the
LQ factorization of the matrix $J$.

We have described algorithm \ref{algo:minmax:family} for a convex
polytope. However, it can be modified for other integration domains as
well. For instance, if $\Omega$ is a solid sphere, then the constraints
result in a family of parabolas and a completely analogous algorithm can
be applied to optimize the parameter $t$.

The constrained LS Newton method is summarized in algorithm~\ref{algo:corrector}.

\begin{algorithm}
\caption{Maps $\tilde \zb$ on a nearby point $\bar \zb$ on $Z_n$.}
\label{algo:corrector}
\begin{algorithmic}
\While {$\Norm{f(\tilde \zb)} > TOL$}
\State Setup $\fb$, $\gb$ and $J$, and compute the LQ-factorization
     of $J$.
\State Calculate $\Delta \zb^f$ and $\Delta \zb^g$ in \eqref{delta:zfg}.
\State Determine $t$ from algorithm~\ref{algo:minmax:family}.
\State Set $\tilde \zb \leftarrow \tilde \zb + \Delta \zb^f + t \Delta \zb^g$ 
\EndWhile
\end{algorithmic}
\end{algorithm}

\section{Node Elimination}\label{sec:nodeelim}
In this section we describe a method to reduce the number
of points in a cubature rule while maintaining its degree. It consists
of determining points on $Z_n$
that intersect with the coordinate planes $w_k = 0$. Obviously, 
if a weight vanishes in \eqref{cubaturerule}, then the corresponding node does not
have to be included in the cubature rule \eqref{cubaturerule}. Thus we
have found another solution of \eqref{def:polysys} with $n-1$ nodes
that has the same degree of precision. The elimination procedure is
then repeated until no further quadrature rules with vanishing nodes
can be found.

We start with a point $\bar \zb$ on $Z$ with positive weights and
use a predictor-corrector type approach to find another point on $Z_n$ where
one of the $w_k$'s  vanishes. The predictor step consists of linearizing $\fb$ at
$\bar \zb$ and then computing the nearest points on the tangent space
that intersect with the one of the hyperplanes $w_k=0$. The $k$-th
node is then eliminated and the point $\tilde \zb$ is mapped to 
$Z_{n-1}$ using algorithm~\ref{algo:corrector}.
The predictor-corrector
method is restarted to eliminate further nodes. This approach
guarantees that all nodes remain in the domain and all weights are non-negative.

We now describe the predictor step in more detail. To obtain the tangent
space of $Z_n$ at $\bar \zb$ consider the linearization of the function
$\fb$ at $\bar \zb$
\begin{equation*}
  \fb(\bar \zb + \Delta \zb) = \fb(\bar \zb) + J \Delta \zb +
  \bigO(\abs{\Delta \zb}^2),
\end{equation*}
where $J = D\fb(\bar \zb)$.
Since the first term on the right hand side vanishes, the tangent space is defined by
$\zb = \bar \zb + \Delta \zb$, where $\Delta \zb \in \Null(J)$.
Since the matrix $J$ is underdetermined, this nullspace is
nontrivial. An orthonormal basis can be found with
the full LQ-factorization of $J$
\begin{equation}\label{lq:factoriz}
  J = \Big[L, 0 \Big] \left[\begin{array}{c}\tilde Q\\ \hat Q \end{array} \right]
\end{equation}
where $L\in \R^{M\times M}$ is lower triangular, $\tilde Q \in
\R^{M\times N}$ and $\hat Q \in\R^{N-M\times N}$. The rows of $\hat Q$
form an orthogonal basis of $\Null(J)$. 

To eliminate the $k$-th node, consider the nearest point $\tilde \zb$ 
in the intersection of the tangent space with the hyperplane
$w_k=0$. If $\tilde \zb = \bar \zb + \Delta \zb$ then $\Delta \zb$
solves the optimization problem
\begin{equation*}
\begin{aligned}
\min \;       &\half \Norm{\Delta \zb}^2 + t\, \gb^T \Delta \zb \\
\Delta \zb:\; & J \Delta \zb =  \mathbf{0}\\
              & w_k + \Delta w_k = 0.
\end{aligned}
\end{equation*}
Here $\gb = \nabla \phi_{\Omega,k}(\bar \zb)$, where $\phi_{\Omega,k}(\cdot)$
is the penalty function that is obtained by excluding the
node $k$ in the summation in equation \eqref{def:penalty:fcn}.

The unknown $\Delta \zb$ can be expressed as a combination of the
orthogonal basis of $\Null(J)$ obtained from the LQ factorization in
\eqref{lq:factoriz}.  Thus there is a vector $\Delta \yb \in
\R^{N-M}$ such that $\Delta \zb = \hat Q^T \Delta \yb$. Substitution
of this vector leads to the optimization problem 
\begin{equation*}
\begin{aligned}
\min \;       &\half \Norm{\Delta \yb}^2 + t\, \hat \gb^T \Delta \yb \\
\Delta \yb:\; & \mb_k^T \Delta \yb = - \bar w_k.
\end{aligned}
\end{equation*}
Here $\mb_k$ is the column of $\hat Q$ corresponding to the $k$-th
weight, and $\hat \gb_k = \hat Q \gb$. It can be seen that the
solution $\Delta \yb_k$ of the above optimization problem is given by 
\begin{equation}\label{def:dyk}
  \Delta \yb_k = \Delta \yb_k^f + t \Delta \yb_k^g,
\end{equation}
where
\begin{equation*}
  \begin{aligned}
    \Delta \yb_k^f &= \mb_k \frac{\bar w_k}{\Norm{\mb_k}^2}\\
    \Delta \yb_k^g &=\left( I - \frac{ \mb_k \mb_k^T}{\Norm{\mb_k}^2}\right) \hat \gb,
    \end{aligned}
\end{equation*}
and thus the predictor is 
\begin{equation*}\label{def:zbzg}
  \tilde \zb_k = \bar \zb + \Delta \zb^f_k + t  \Delta \zb^g_k,
\end{equation*}
where $\Delta \zb^f = \hat Q^T \Delta \yb_k^f$ and 
$\Delta \zb^g = \hat Q^T\Delta \yb_k^g$.
The parameter $t$ is selected to ensure maximal distance from the domain boundary.
If the domain $\Omega$ is a polytope, then this can be achieved with
algorithm~\ref{algo:minmax:family}. 

We compute in a list predictors $\tilde \zb_k$ for all weights $k\in
\{ 1,\dots,n\}$. The predictors that satisfy the constraints and are
close to $\bar \zb$ will be mapped on $Z_{n-1}$. The next
quadrature is the solution that has the greatest distance from the boundary.
The node elimination procedure is summarized in algorithm~\ref{algo:nodeelim}.

\begin{algorithm}
\caption{Node elimination algorithm.}
\label{algo:nodeelim}
\begin{algorithmic}
\State Find a suitable initial cubature rule $\bar \zb$.
\While {$N > M$}
\State Setup $J$, and compute the LQ-factorization.
\For{k=1:n}
\State Compute $\Delta \zb_k^f$ and $\Delta \zb_k^g$ in \eqref{def:zbzg}.
\State Compute $t$ using using algorithm~\ref{algo:minmax:family}.
\State Set $\Delta \zb_k = \Delta \zb^f_k + t  \Delta \zb^g_k$ 
     and $\tilde \zb_k = \bar \zb_k + \Delta \zb_k$.
\EndFor
\State Sort the $\Delta \zb_k$'s that satisfy the constraints such that 
\State $\qquad \Norm{\Delta \zb_1}\leq \Norm{\Delta \zb_2} \leq \dots$.
\State \For{k=1:K}
\State Eliminate the $k$-th node from $\tilde \zb_k$.
\State Use algorithm~\ref{algo:corrector} to map $\tilde \zb_k$ to
$\bar \zb_k \in Z_{n-1}$ 
\EndFor
\State \textbf{Stop} if no $\bar \zb_k \in Z_{n-1}$ could be found.
\State $\bar \zb \leftarrow \bar \zb_k$, where $\bar \zb_k$'s nodes have the greatest distance from
the boundary.
\State $n \leftarrow n-1$
\EndWhile
\end{algorithmic}
\end{algorithm}

The main cost in this algorithm is to set up Jacobians and compute
their QR factorizations.

\section{Strategies for the Initial Cubatures}\label{sec:initialguess}
The discussion so far assumed that an initial cubature rule is known
which may be suboptimal (i.e., it has more nodes than necessary), but
has the desired degree. This section describes how such an initial
cubature can be obtained. Here the goal is to keep the number of nodes
low such that fewer elimination steps have to be taken and the cost of
the linear algebra in the initial stages of the algorithm is reduced.

The scheme is based on the fact that cubature rules for Cartesian product
domains $\Omega = \Omega_1 \times \Omega_2$
can
be obtained by forming tensor products of cubature rules of $\Omega_1$
and $\Omega_2$. Specifically, if $\{\xb_n,w_n\}$ and $\{\yb_m,v_m\}$ are
cubature rules for $\Omega_1$ and $\Omega_2$ that are exact for
$\Pol_p^{d_1}$ and $\Pol^{d_2}_p$, respectively, then
$\{(\xb_n,\yb_m),w_nv_m\}$ is a cubature rule that is exact for
$\Pol^{d_1}_p \times \Pol^{d_2}_p$. Since this polynomial space is
larger than $\Pol^{d_1+d_2}_p$, the node elimination can be performed
with the tensor product rule as an initial guess.

We start the discussion with the case that $\Omega$ is the $d$-dimensional unit cube
\begin{equation*}
  C_d = \big\{ (x_1,\dots,x_d) : 0\leq x_i \leq 1, i=0,\dots,d \big\}.
\end{equation*}
If one were to use tensor products of degree-$p$ Gauss Legendre rules
as the initial quadrature,
one would obtain $(\floor{p/2}+1)^{d}$ nodes. The
rapid growth of this number makes this strategy unfeasible even for
moderate values of $p$ and $d$.

Instead, our implementation uses an iterative scheme to obtain the
initial guess for the domains of interest, thereby significantly
decreasing the initial number of nodes. We build the rule
incrementally, by starting with $C_1\times C_1$, and running node
elimination. The resulting rule is then tensored with the $C_1$-rule to obtain an
initial guess for $C_3$. The scheme is repeated until 
the desired dimension is reached. With this approach the number of
nodes of the initial guess in the last step is greatly reduced over a
$d$-fold tensor product of rules for $C_1$.

We now turn to the d-dimensional simplex
\begin{equation*}
  T_d = \big\{ (x_1,\dots,x_d) : 0\leq x_d \leq \dots \leq x_1 \leq 1 \big\}
\end{equation*}
and to Cartesian products of cubes and simplices. 

The well known Duffy transformation can be used to transform a cube
into a simplex. For the $d+1$-dimensional simplex it can be defined
recursively as follows
\begin{equation}\label{duffy}
\begin{aligned}
x_1 &=\xi, \\ 
 x_2 &=\xi   \eta_1,  \\ 
  x_3 &=\xi \eta_2, \\
 \vdots \\
  x_{d+1} &=\xi \eta_d.
\end{aligned}
\end{equation}
Here $\xi \in [0,1] = C_1$ and $\etab = (\eta_1,\dots,\eta_d)\in
T_{d}$. Thus \eqref{duffy} transforms $C_1 \times T_d$ into $T_{d+1}$ 
with Jacobian $J = \xi^{d}$. Now $T_d$ can be transformed with another
Duffy transform and repeating this leads to
a transformation from $C_{d+1}$ to $T_{d+1}$. However, for the
construction of initial quadrature rules it is more convenient to work with
one transformation at a time.

For instance, for the triangle (i.e. $d=1$), we use Gauss-Jacobi rules
for the integral over the $(\xi,\eta)$ variables. 
Specifically, if
\begin{equation*}
  \int_0^1 f(\eta) d\eta \approx \sum_{\ell=1}^q f(y_\ell) v_\ell
  \quad\mbox{and}\quad
\int_0^1 f(\xi) \xi d\xi \approx \sum_{k=1}^q f(x_k) w_\ell
\end{equation*}
are the quadrature rules of degree $p=2q-1$ for weight function $w(\eta)=1$ and
$w(\xi)=\xi$, respectively, 
then an integral over $T_2$ can be approximated as follows 
\begin{equation}\label{quad:S2:Duffy}
  \int_{T_2} \varphi(\xb) d\xb = \int_{0}^{1} \int_{0}^{1}
  \varphi(\xi,\xi\eta) \xi d\eta d\xi
  \approx \sum_{k=1\atop \ell=1}^p \varphi(x_k,x_ky_\ell) w_k v_\ell.
\end{equation}
Note that the Jacobian $\xi$ determines the choice of rule and affects
the nodes $x_k$ and weights $w_k$.  If $\varphi \in \Pol_p^2$, then
$(\xi,\eta) \mapsto \varphi(\xi, \xi \eta)$ is a polynomial in
$\Pol_p^1\times \Pol_p^1$ and thus the quadrature rule
\eqref{quad:S2:Duffy} is exact.  This tensor product rule can now be
used as the initial rule in the node-elimination procedure.

Similar to the cube, the 
$d+1$-dimensional simplex rule is built by recursion using previously
generated rules. 
Thus, if $\{\xb_k,w_k\}$ is a degree-$p$ Gauss-Jacobi rule for
$\int_0^1 f(\xi)\xi^d d\xi$  
and $\{\yb_l,v_l \}$ is a degree-$p$ rule for $T_d$, then
\begin{equation*}
  \int_{T_{d+1}} \varphi(\xb) d\xb
= \int_0^1 \int_{T_d} \varphi\big(\xi, \xi \etab \big) \xi^d d\eta  d\xi 
  \approx \sum_{k,l} \varphi\big(x_k, x_k \yb_l \big) w_k v_l
\end{equation*}
is a degree-$p$ rule for $T_{d+1}$.

For tensor products of cubes and simplices the analogous procedures
can be applied, where some care must be a taken to use the appropriate
Gauss-Jacobi rule to compensate for the Jacobian of the Duffy
transformation.

\begin{figure}[hbt]
\begin{center}

\begin{tikzpicture}[        
ellipsenode/.style={
	ellipse, 
	fill=orange,
	inner sep=4pt,
	text width=5em, 
	text=white,
	text centered
},
rectnode/.style={
	sharp corners, 
	fill=white,
	inner sep=4pt,
	text=black,
	text width=2em, 
	text centered
},
graynode/.style={
	draw=gray,
	very thick,
	inner sep=4pt,
	text=gray,
	text centered
},
myline/.style={
	draw=#1,
	line width=1.5pt,
},
terminator/.style = {shape=rounded rectangle
},  
]
	\node[draw,terminator, cyan] (R1C1) at (1.5,-2) {$\mathbf{I}$};
	\node[draw,rectnode] (R1C2) at (3,-2) {$C_2$};
	\node[draw,rectnode] (R1C3) at (5,-2) {$\mathbf{C_2}$};
	\node[draw,rectnode] (R1C4) at (6.5,-2) {$C_3$};
	\node[draw,rectnode] (R1C5) at (8.5,-2)
	 {$\mathbf{C_3}$};
	\node[draw,rectnode] (R1C6) at (10,-2) {$C_4$};
	\node[draw,terminator, purple] (R1C7) at (12,-2)
	{$\mathbf{C_4}$};

	\node[draw,rectnode] (R2C2) at (3,-4) {$T_2$};
	\node[draw,rectnode] (R2C3) at (5,-4) {$\mathbf{T_2}$};
	\node[draw,rectnode] (R2C4) at (6.5,-4) {$T_2C_1$};
	\node[draw,rectnode] (R2C5) at (8.5,-4)
	{$\mathbf{T_2C_1}$};
	\node[draw,rectnode] (R2C6) at (10,-4) {$T_2 C_2$};
	\node[draw,terminator, purple] (R2C7) at (12,-4)
	{$\mathbf{T_2C_2}$};

	\node[draw,rectnode] (R3C4) at (6.5,-6) {$T_3$};
	\node[draw,rectnode] (R3C5) at (8.5,-6)
	{$\mathbf{T_3}$};
	\node[draw,rectnode] (R3C6) at (10,-6) {$T_3C_1$};
	\node[draw,terminator, purple] (R3C7) at (12,-6)
	{$\mathbf{T_3C_1}$};
	
    \node[draw,rectnode] (R4C6) at (10,-8) {$T_2T_2$};
	\node[draw, terminator, purple] (R4C7) at (12,-8) {$T_2 T_2$};

	\path [line] (R1C2)--node[text width=0.5cm,above] {NE} (R1C3);
	\path [line] (R1C4)--node[text width=0.5cm,above] {NE} (R1C5);
	\path [line] (R1C6)--node[text width=0.5cm,above] {NE} (R1C7);
	
	\path[every node/.style={sloped,anchor=south,auto=false}]
	(R1C1) edge[->, double distance=1pt, >=latex', bend left=45] node{$\otimes I$} (R1C2)
	
	(R1C3) edge[->, double distance=1pt, >=latex', bend left=45] node{$\otimes I$}(R1C4)
	(R1C5) edge[->, double distance=1pt, >=latex', bend left=45] node{$\otimes I$} (R1C6);

	\path [line] (R2C2)--node[text width=0.5cm,above] {NE} (R2C3);
	\path [line] (R2C4)--node[text width=0.5cm,above] {NE} (R2C5);
	\path [line] (R2C6)--node[text width=0.5cm,above] {NE} (R2C7);
	

	 \draw[->=gray, solid] (R2C3) node[text height=3.5cm,text width=0.5,below]{$\otimes T_2$} |- (R4C6) 
	;
	
	\path[every node/.style={sloped,anchor=south,auto=false}]
	(R2C3) edge[->, double distance=1pt, >=latex', bend left=45] node{$\otimes I$} (R2C4)
	(R2C5) edge[->, double distance=1pt, >=latex', bend left=45] node{$\otimes I$} (R2C6);

	\path [line] (R3C4)--node[text width=0.5cm,above] {NE} (R3C5);
	\path [line] (R3C6)--node[text width=0.5cm,above] {NE} (R3C7);
	
	\path [line] (R4C6)--node[text width=0.5cm,above] {NE} (R4C7);

		\path[every node/.style={sloped,anchor=south,auto=false}]
	(R3C5) edge[->, double distance=1pt, >=latex', bend left=45] node{$\otimes I$} (R3C6);
	
	\path[every node/.style={sloped,anchor=south,auto=false}]
	(R1C2) edge      node{Duffy} (R2C2) ;         
	\path[every node/.style={sloped,anchor=south,auto=false}]
	(R2C4) edge      node{Duffy} (R3C4) ;

\end{tikzpicture}
\end{center}
\caption{Construction scheme for quadrature rules for the four
  dimensional polytopes $C_4$, $T_2\times C_2$, $T_3\times C_1$ $T_2\times T_2$}
\label{fig:multiDscheme}
\end{figure}
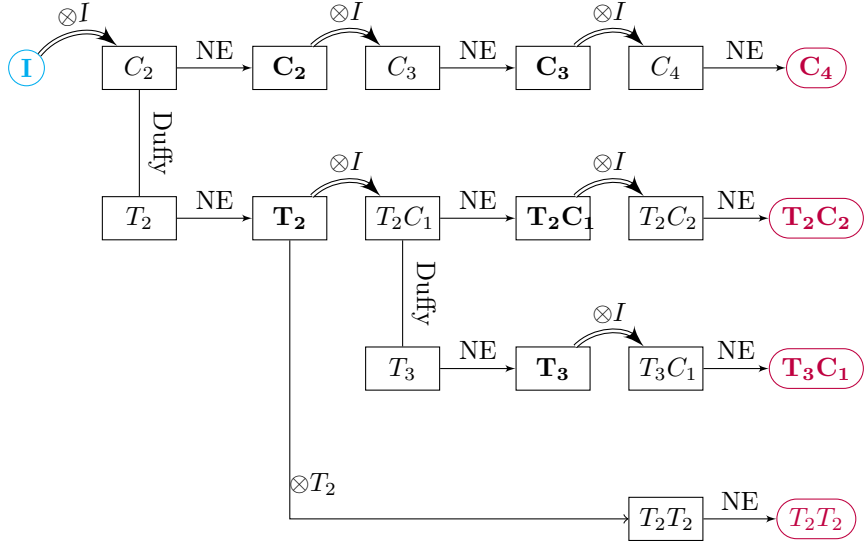

Figure~\ref{fig:multiDscheme} illustrates the procedure for various
polytopes in four dimensions.  It begins with Gaussian cubature on the
interval (shown in blue). Then the tensor product of Gaussian
quadrature on $I \otimes I$ is applied, which is an initial guess for
$C_2$. Subsequently, depending on the domain of interest, we either
proceed by running the Node Elimination algorithm, or applying the Duffy
transformation to obtain $T_2$. Figure \ref{fig:multiDscheme} describes how to
derive initial guesses for $C_4, \ T_2\times C_2,\ T_3\times C_1$ and
$T_2\times T_2$, before performing the final node elimination step,
which is shown in red.


\section{Numerical Examples}\label{sec:numresult}
We have implemented the method in C++ using double precision
arithmetic and the LAPACKE interface to the multithreaded
implementation of LAPACK library for the dense linear algebra. The
tests were run on a single process of a standard desktop system with
an Intel I7 processor and 32 gigabyte of memory. All shown examples
could be run with much less memory usage than our system had.

For the basis functions in
\eqref{def:vandermonde} we use tensor products of Legendre polynomials
for the cubes and the orthogonal triangular polynomials in
\cite{koornwinder75}. The construction in this paper can be
extended to obtain orthogonal polynomials for simplices of any
dimension, which we implemented in our code.

The number of iterations in the corrector step is always
small. Throughout most of the node elimination algorithm, the
corrector converges after one to three iterations. Only in the final
stages, when the system is less overdetermined, the typical number of
iterations increases to three to seven, in rare cases up to ten. The
elimination procedure is stopped if either 
the optimal number of nodes has been reached or  
if algorithm~\ref{algo:minmax:family} is unable to find a
value of $t$ for all predictors. The tolerance in 
algorithm~\ref{algo:corrector} is $TOL=10^{-14}$.

Table~\ref{tab:twoDrules} displays data about the final cubature
rules obtained with our node elimination algorithm
for the domain $T_2$. Here $\ntensor$ is the number of points in the
tensor product rule \eqref{quad:S2:Duffy}, which is the initial
cubature rule,  and $\nelim$ is the number of
points in the final rule.
The optimal number of nodes and the efficiency ratio is given by
\begin{equation*}
  \nopt= \left\lceil \frac{\mbox{dim }\Pol^d_p}{(d+1)} \right\rceil
  \quad\mbox{and}\quad \iopt = \frac{\nopt}{\nelim}.
\end{equation*} 
The table contains only odd degree rules, because there are no even
degree tensor product rules. Even degree rules could have been easily obtained
by starting with the next degree odd tensor product rule, and
enforcing only the even degree in the moment equations.

\begin{table}
\begin{center}
\begin{tabular}{lcccccccccccccc}
degree        & 3   & 5    & 7    & 9    & 11   & 13   & 15  \\ 
\hline
$\ntensor$    & 4   & 9    & 16   & 25   & 36   & 49   & 64  \\ 
$\nelim$      & 4   & 7    & 12   & 19   & 27   & 38   & 47  \\ 
$\nopt$       & 4   & 7    & 12   & 19   & 26   & 35   & 46  \\ 
$\iopt$       &1.00 & 1.00 & 1.00 & 1.00 & 0.96 & 0.92 & 0.98 \\
\\
degree        & 17  & 19   & 21   & 23   & 25   & 27   & 29 \\       
\hline
$\ntensor$    & 81  & 100  & 121  & 144  & 169  & 196  & 225 \\    
$\nelim$      & 58  & 74   & 86   & 102  & 119  & 142  & 163\\    
$\nopt$       & 57  & 70   & 85   & 100  & 117  & 136  & 155 \\      
$\iopt$       &0.98 &0.94  &0.99  &0.98  & 0.98 & 0.96 & 0.95\\
\end{tabular}
\caption{Data for the cubature rules on the triangle}
\label{tab:twoDrules}
\end{center}
\end{table}

One can see that our algorithm found rules that are very close to 
optimal. One should note that rules for $T_2$ of similar quality have
been obtained by other authors, including \cite{xiao-gimbutas10}. 

\begin{figure}[hbt]
\begin{center}
\includegraphics[width=2.35in]{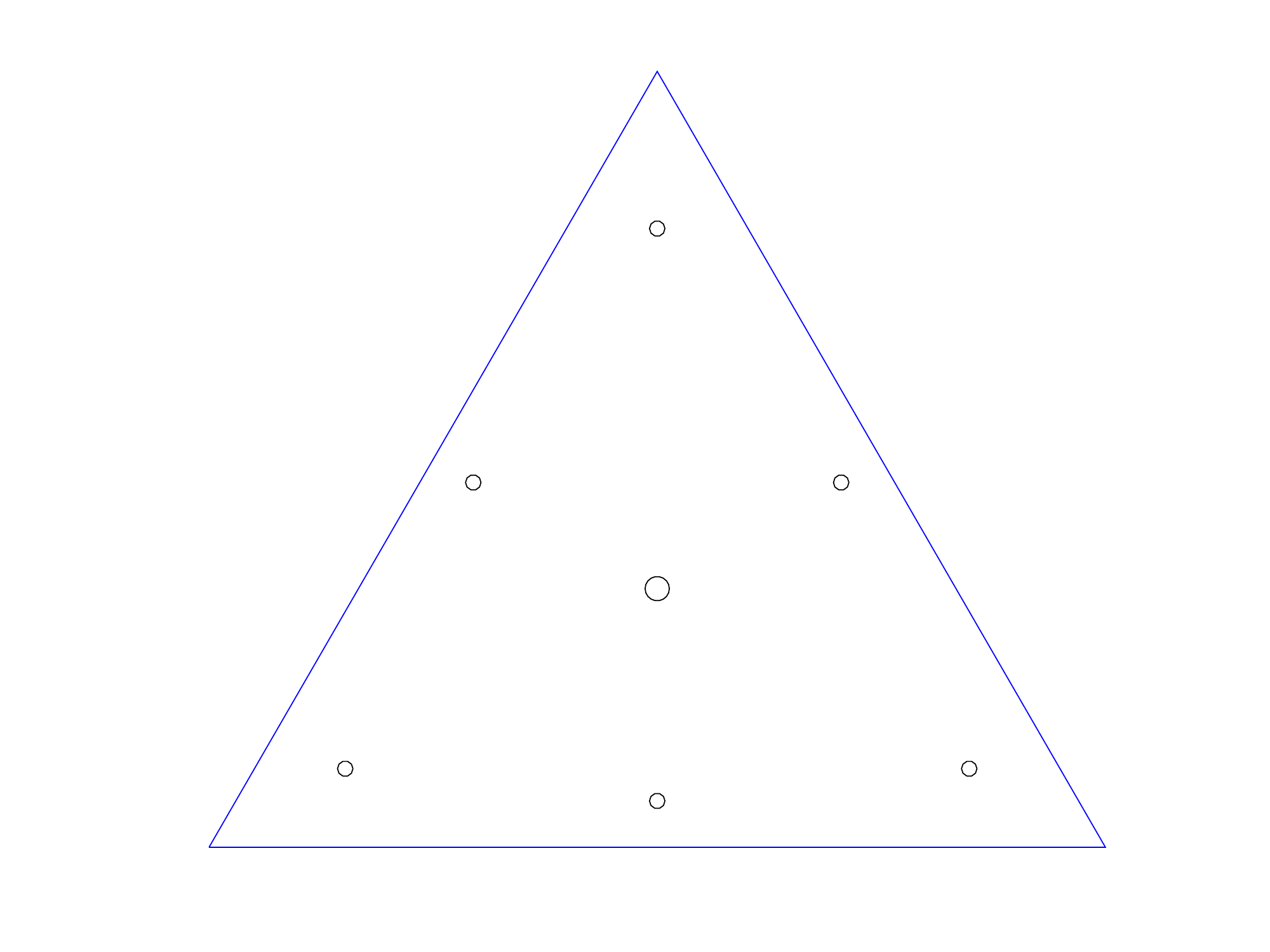}
\includegraphics[width=2.35in]{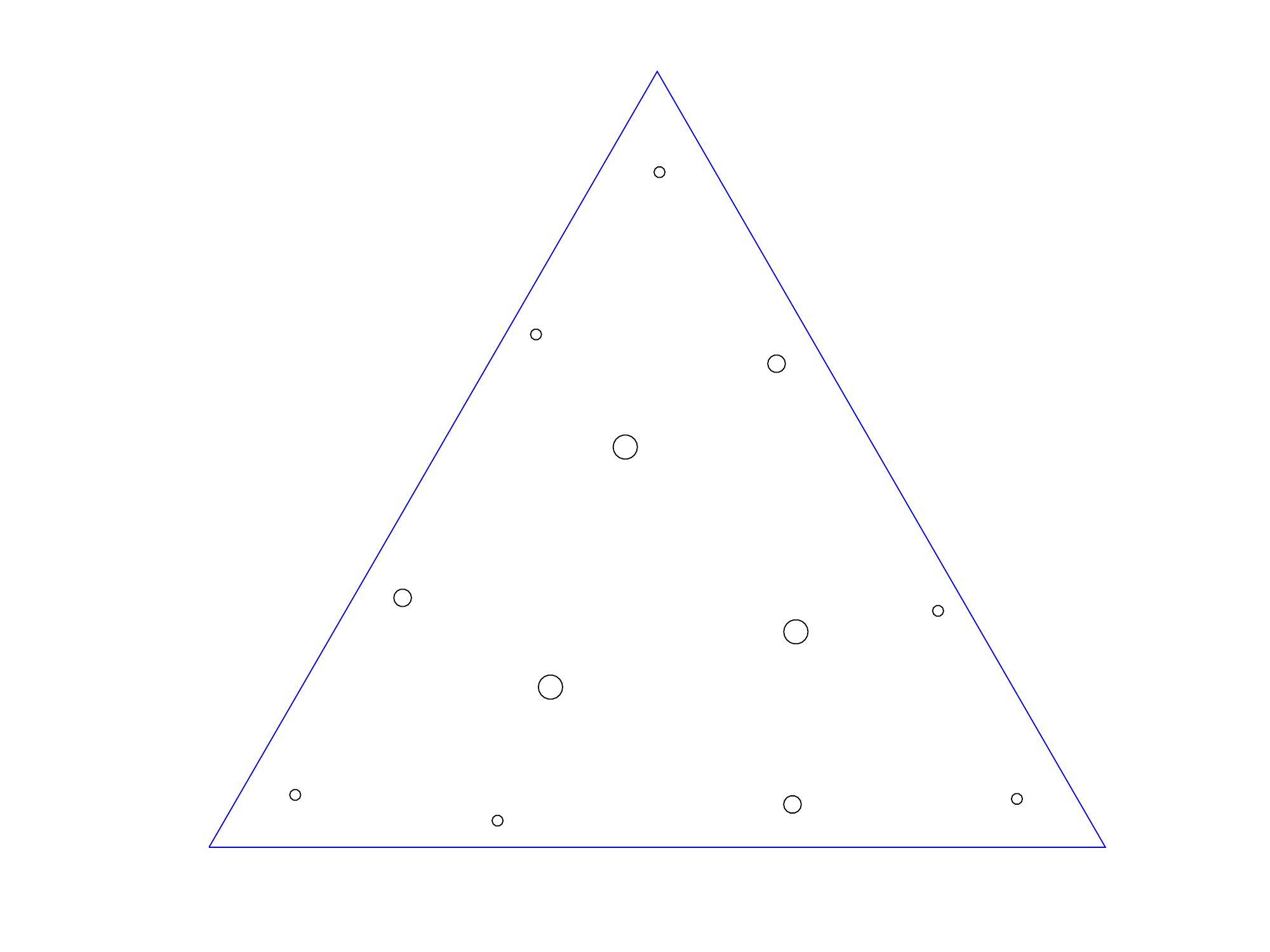}\\
\includegraphics[width=2.35in]{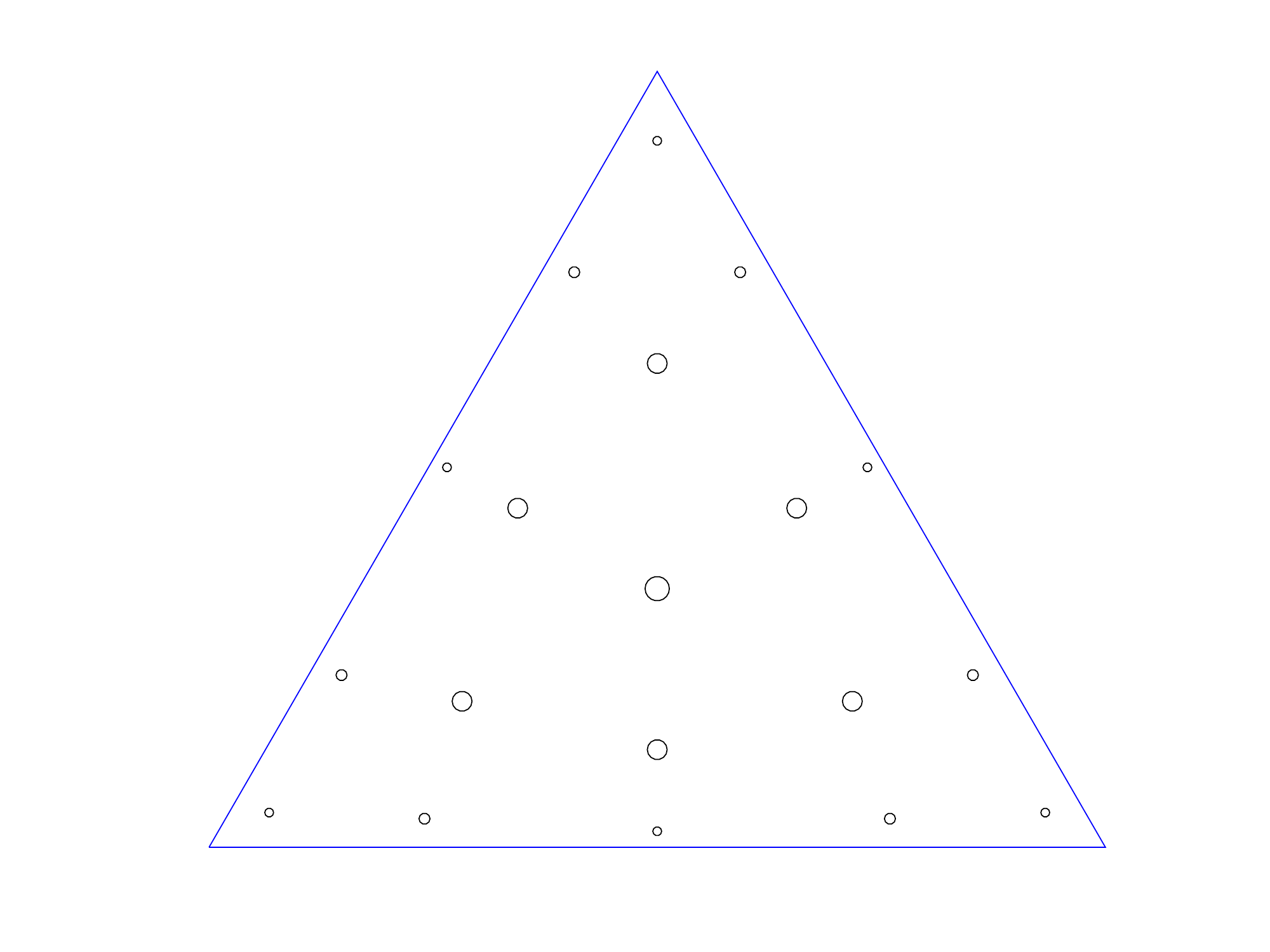}
\includegraphics[width=2.35in]{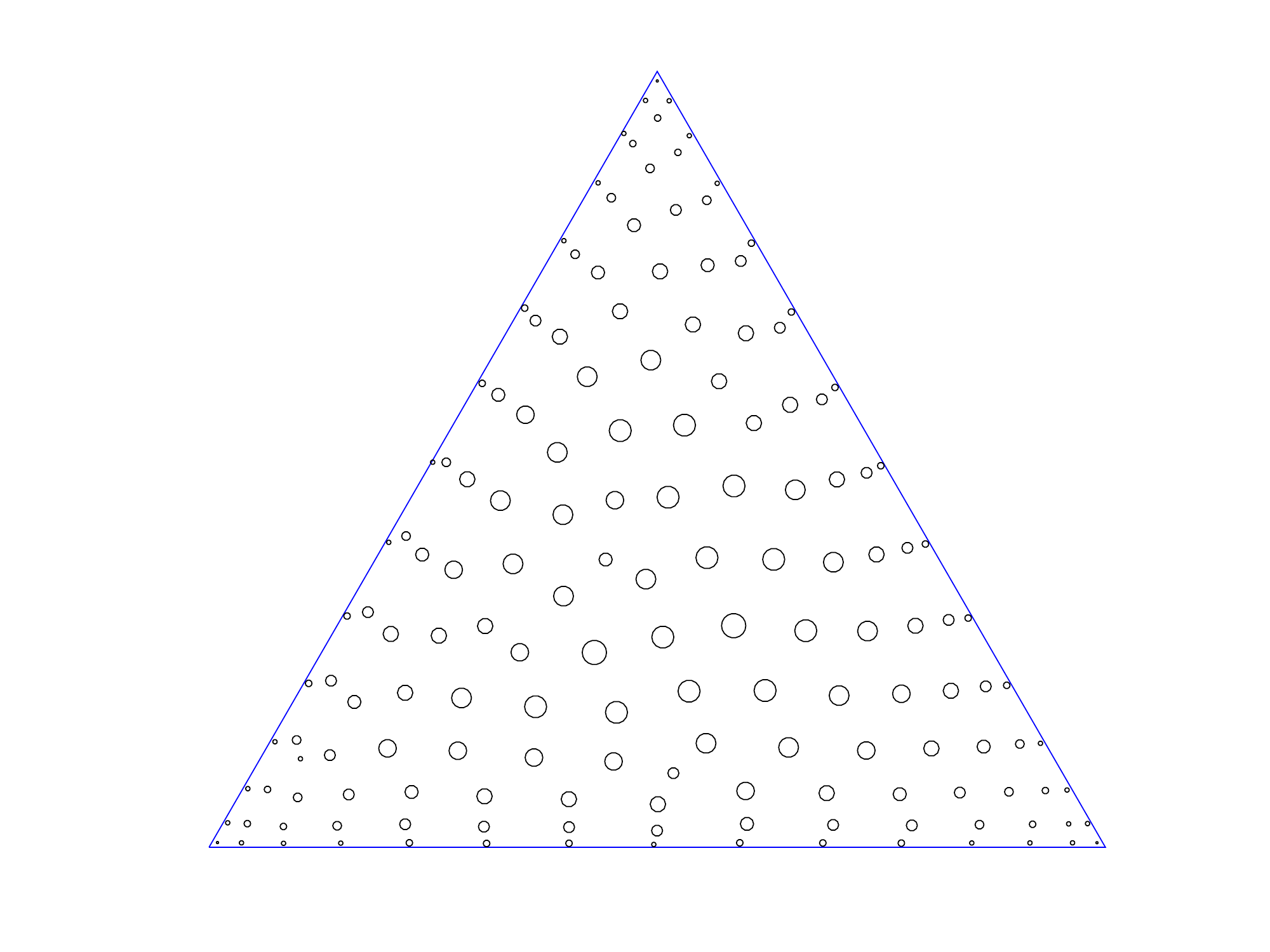}
\end{center}
\caption{Quadrature nodes for a triangle. The degrees are 5, 7, 9 and 29.}
\label{fig:triangles}
\end{figure}

Our implementation does not exploit or enforce any symmetry, therefore
the computed rules are in general not symmetric.
Figure~\ref{fig:triangles} displays for some examples the location of
the nodes and size of the weights, indicated by the radius of the
circles. For the figure, we have mapped $T_2$ to an equilateral
triangle to better assess the symmetries of the nodes. The rule of
degree five is fully symmetric and is identical to a rule that has
been found by analytical methods, see~\cite{stroud71}. It is also the
only rule we found where $M=N$ in the moment equations. This implies
that the solution is an isolated point, and it makes sense to
calculate the difference of the analytic solution with the numerically
obtained one. We have found that the maximal difference is
$1.14 \cdot 10^{-16}$, which is well within machine precision.

The rule of degree seven only has the three rotational symmetries, but
none of the reflectional symmetries of the equilateral triangle. The
rule of degree nine is symmetric up to seven digits of accuracy. The
nearby symmetric version is known, see
Dunavant \cite{dunavant85b}. Both rules
solve the moment equations to machine precision. Note that while this
rule is optimal, the moment equations are still underdetermined,
because $N=3\cdot 19 = 57$ and $M= \half \cdot 10\cdot 11 = 55$. Thus
our rule and Dunavant's rule are nearby solutions on a two dimensional
manifold.

The higher degree rules in table~\ref{tab:twoDrules} are not
symmetric, but the nodes are generally following a symmetric
distribution and the weights tend to be larger in the interior of the
domain. This happens despite the fact that the nodes of the initial
rules are more concentrated near the origin.

For the lower degree rules it can be determined analytically whether
and how many fully symmetric rules exist,
see~\cite{papanicolopulos15}. Thus it is interesting to compare the
theoretically predicted PI rules with those we obtained in
table~\ref{tab:twoDrules}. In addition to the already mentioned degree
five and nine rules, the paper~\cite{papanicolopulos15} states that
there are four
fully symmetric PI rules of degree 11 with 30 nodes, and two rules of
degree 13 with 37 points and one rule of degree 15. Our degree 11 rule
has fewer nodes but is not symmetric, while the degree 13 rule has one
more node and the degree 15 rule has the same number of nodes,
albeit no symmetries. There is no symmetric PI rule of degree 7, but
our algorithm was able to find an optimal rule with partial symmetry.

\begin{table}
\begin{center}
\begin{tabular}{lccccccccccccccc}
degree      & 5     & 7    & 9   & 11  & 13   & 15  & 17  & 19   &21 & 23\\ 
\hline                                                                
$\ntensor$  & 27    & 64   & 125 & 216 & 343  & 512 & 729 & 1000 &1331 &1728\\ 
$\nelim$    & 14    & 31   & 57  & 95  & 143  & 206 & 288 & 390  &510  &653\\ 
$\nopt$     & 14    & 30   & 55  & 91  & 140  & 204 & 285 & 385  &506  & 650\\ 
$\iopt$     & 1.00  & 0.97 &0.96 &0.96 &0.98  &0.99 &0.99 &0.99  &0.99 &0.99 \\
\end{tabular}
\caption{Data for the cubature rules on the 3D simplex}
\label{tab:S3rules}
\end{center}
\end{table}
The results for the domain $T_3$ are shown in table
\ref{tab:S3rules}. The encyclopedia of cubature rules~\cite{cools03}
lists a PI rule of degree 6 with 29 nodes, and higher degree rules
with indefinite weights.  The
paper~\cite{jaskowiec-sukumar19} reports a numerically generated rule
of degree 19 with 392 and a rule of degree 20 with 448 nodes.

\begin{table}
\begin{center}
\begin{tabular}{lccccccc}
degree      & 3     & 5     & 7    & 9   & 11  & 13   & 15    \\ 
\hline                                                                
$\ntensor$    &  16   & 81   & 256 & 625 & 1296 & 2401 & 4096  \\ 
$\nelim$    &  8    & 26   & 68  & 150 & 283  & 497  &  787  \\ 
$\nopt$    &  7    & 26   & 66  & 143 & 273  & 476  &  776  \\ 
$\iopt$  & 0.88  & 1.00 & 0.97& 0.95& 0.96 & 0.96 &  0.99  \\
\end{tabular}
\caption{Data for the cubature rules for $T_4$.}
\label{tab:S4rules}
\end{center}
\end{table}

\begin{table}
\begin{center}
\begin{tabular}{lccccccc}  
degree      & 3     & 5     & 7    & 9   & 11  & 13   & 15   \\ 
\hline                                                                
$\ntensor$    &  16   & 81   & 256 & 625 & 1296 & 2401 & 4096 \\ 
$\nelim$    &  8    & 26   & 67  & 146 & 277  & 478  &  781 \\ 
$\nopt$    &  7    & 26   & 66  & 143 & 273  & 476  &  776 \\ 
$\iopt$  & 0.88  & 1.00 & 0.99& 0.98& 0.99 & 0.99 & 0.99 \\
\end{tabular}
\caption{Data for the cubature rules for $T_2\times T_2$.}
\label{tab:S2S2rules}
\end{center}
\end{table}

Finally, tables~\ref{tab:S4rules} and~\ref{tab:S2S2rules} show the data
of the numerically computed rules for the four dimensional domains
$T_4$ and $T_2\times T_2$. For comparison, the highest PI
$T_4$-quadrature rule in \cite{cools03} has degree five and 31 nodes.

The quadrature rules can be found in the git-hub repository \cite{slobodkins}.

\subsection*{Quadrature errors} To test the cubature rules over the
simplices we consider the integral 
\begin{equation}\label{testSd}
I = \int\limits_{\Omega} \frac{\exp(\mathbf{a}^T\xb) -1}{\mathbf{a}^T\xb}  d\xb 
\end{equation}
where $\Omega \in \{ T_2, T_3, T_4, T_2\times T_2\}$. Note that the
integrand is an analytic function in all variables, even if
$\mathbf{a}^T\xb$ vanishes. The vector $\mathbf{a}$ is chosen randomly
from the hypercube $[-25,25]^d/\sqrt{d}$, where $d$ is the dimension of
$\Omega$. To compute the reference value we express the integrand in form
of an additional integral, which is computed numerically. Thus
\begin{equation}\label{testSd:int:t}
  I = \int_0^1 \int\limits_{\Omega} \exp(t \mathbf{a}^T\xb) d\xb\, dt
  \approx \sum_k \int\limits_{\Omega} \exp(t_k \mathbf{a}^T\xb) d\xb\, w_k,
\end{equation}
where $t_k,w_k$ are the Gauss-Legendre points. The integral in
\eqref{testSd:int:t} over the
exponential function can be determined analytically for the considered
domains and we use a
quadrature rule for the $t$-variable of sufficient degree to be sure that the integral has
been determined to machine precision. We compute one thousand samples
of the vector $\mathbf{a}$ and report the maximal errors.
Comparisons of two different quadrature rules always involve the
same $\mathbf{a}$-vectors.

\begin{figure}[hbt]
\begin{center}
\includegraphics[width=2.2in]{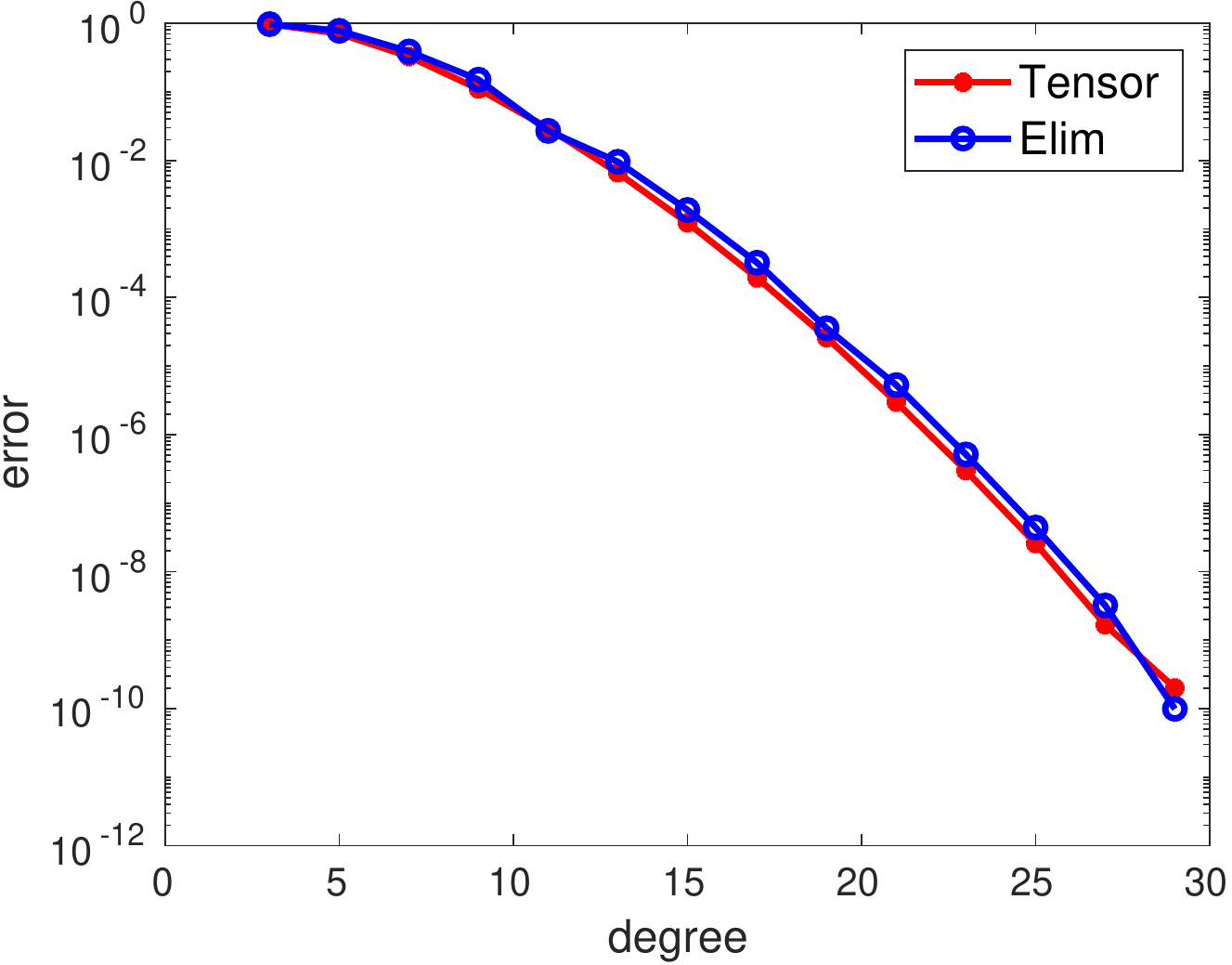}
\includegraphics[width=2.2in]{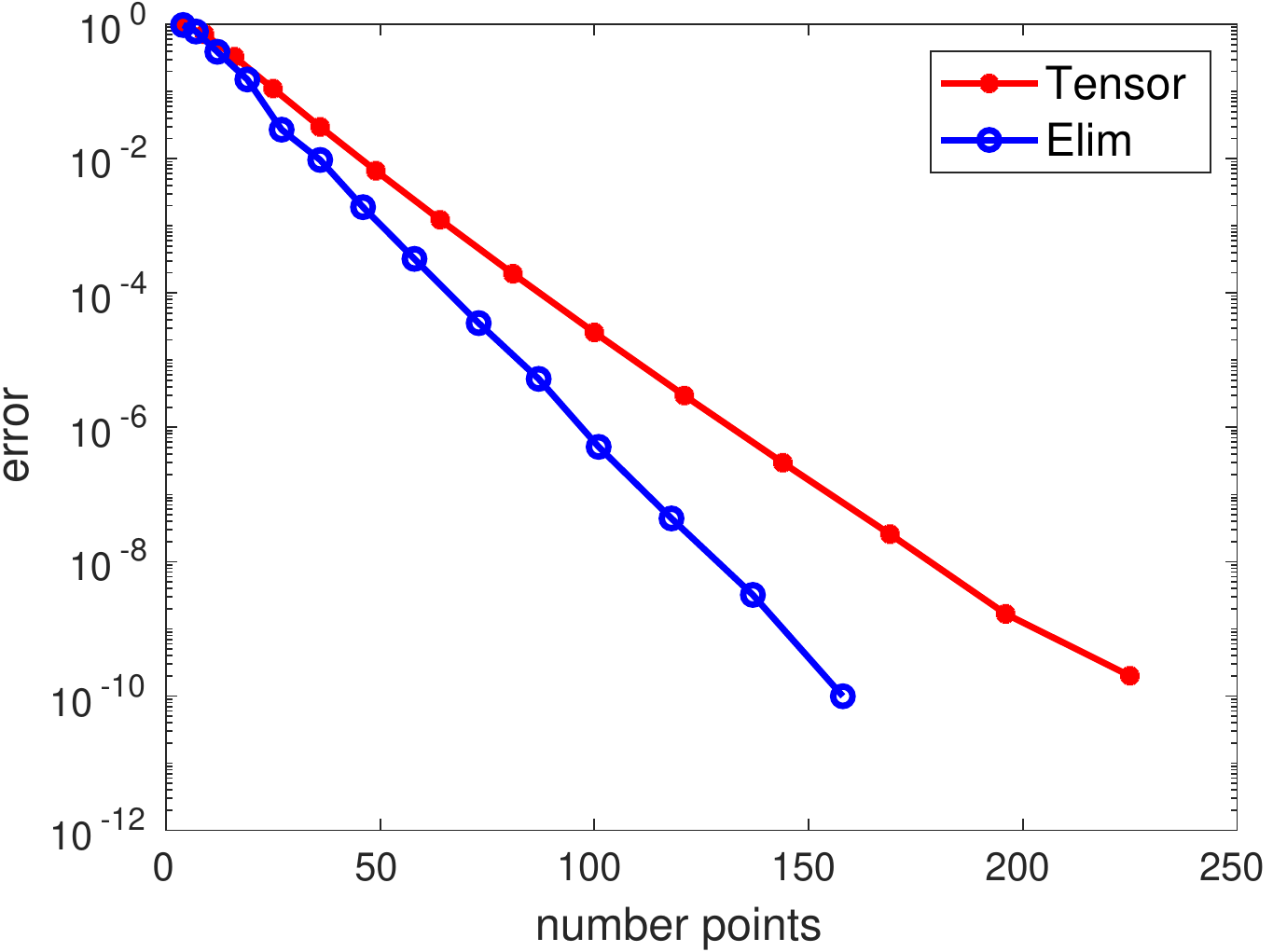}
\end{center}
\caption{Comparison of the tensor product rules with
  the rules of table~\ref{tab:twoDrules}. 
  $\Omega=T_2$. Left: error vs. degree. Right: error vs. number points. }
\label{fig:twoD}
\end{figure}

When $\Omega=T_2$  
the standard approach is to approximate the transformed integral 
using the rule \eqref{quad:S2:Duffy}. Here
we compare this rule with the cubature rules of table~\ref{tab:twoDrules}. 
This is equivalent to comparing the initial quadrature rule
with the final rule in the node elimination process. The results are
shown in 
figure \ref{fig:twoD}, where the left plot shows the error versus the
degree, and the right plot shows the error versus the number of nodes.

\begin{figure}[hbt]
\begin{center}
\includegraphics[width=2.2in]{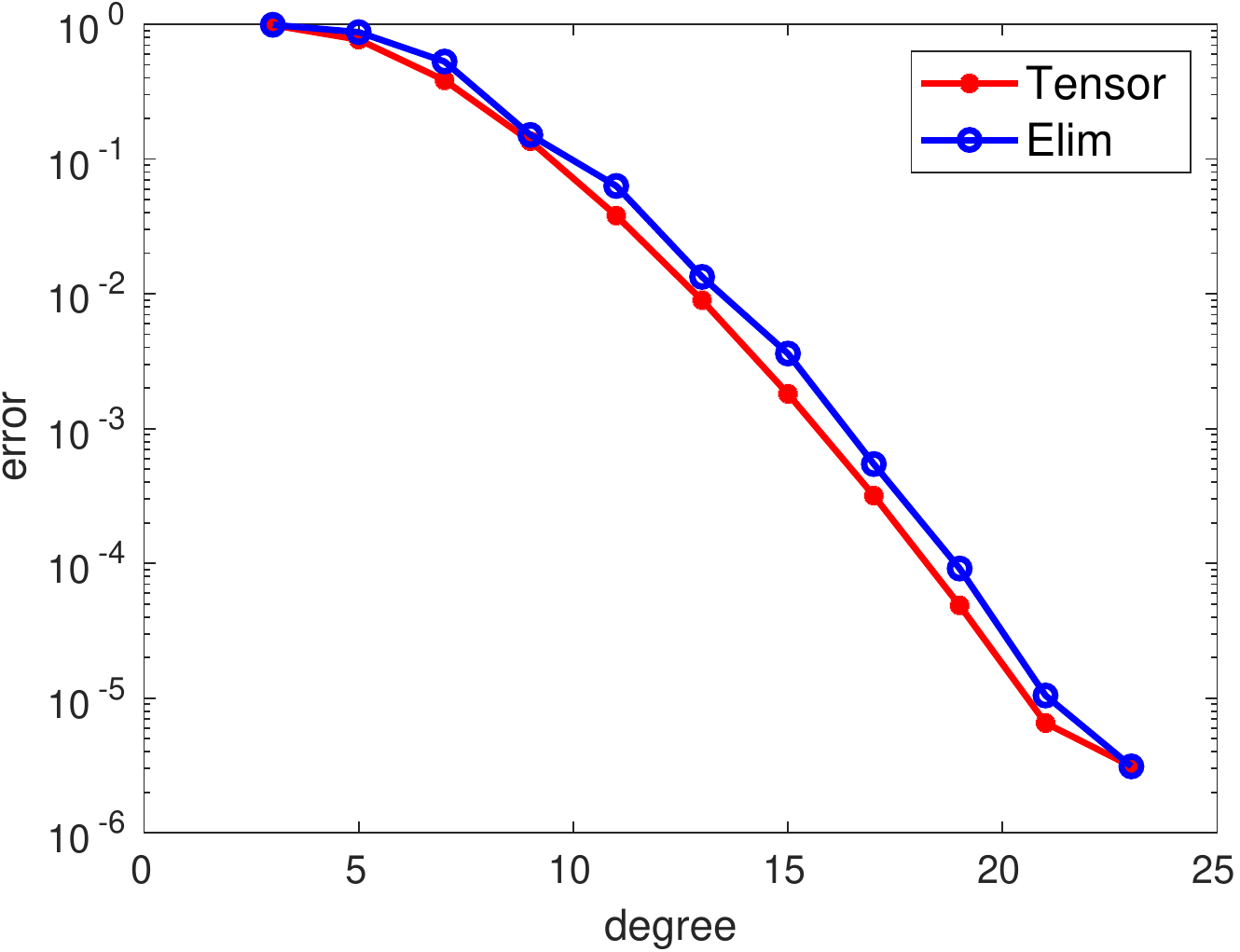}
\includegraphics[width=2.2in]{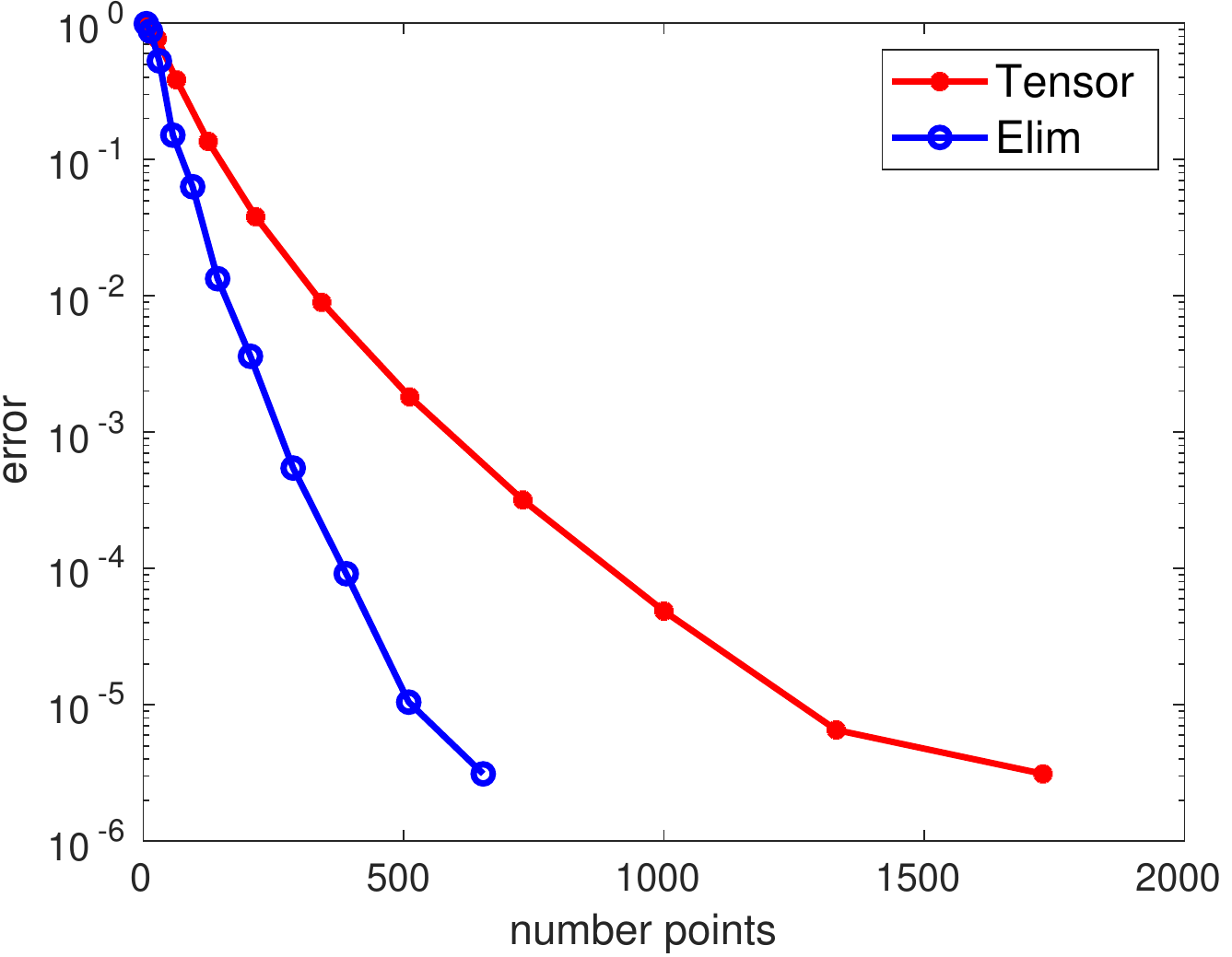}
\end{center}
\caption{Comparison of the tensor product rules with
  the rules of table~\ref{tab:S3rules}. 
  $\Omega=T_3$. Left: error vs. degree. Right: error vs. number points.}
\label{fig:threeD}
\end{figure}

\begin{figure}[hbt]
\begin{center}
\includegraphics[width=2.2in]{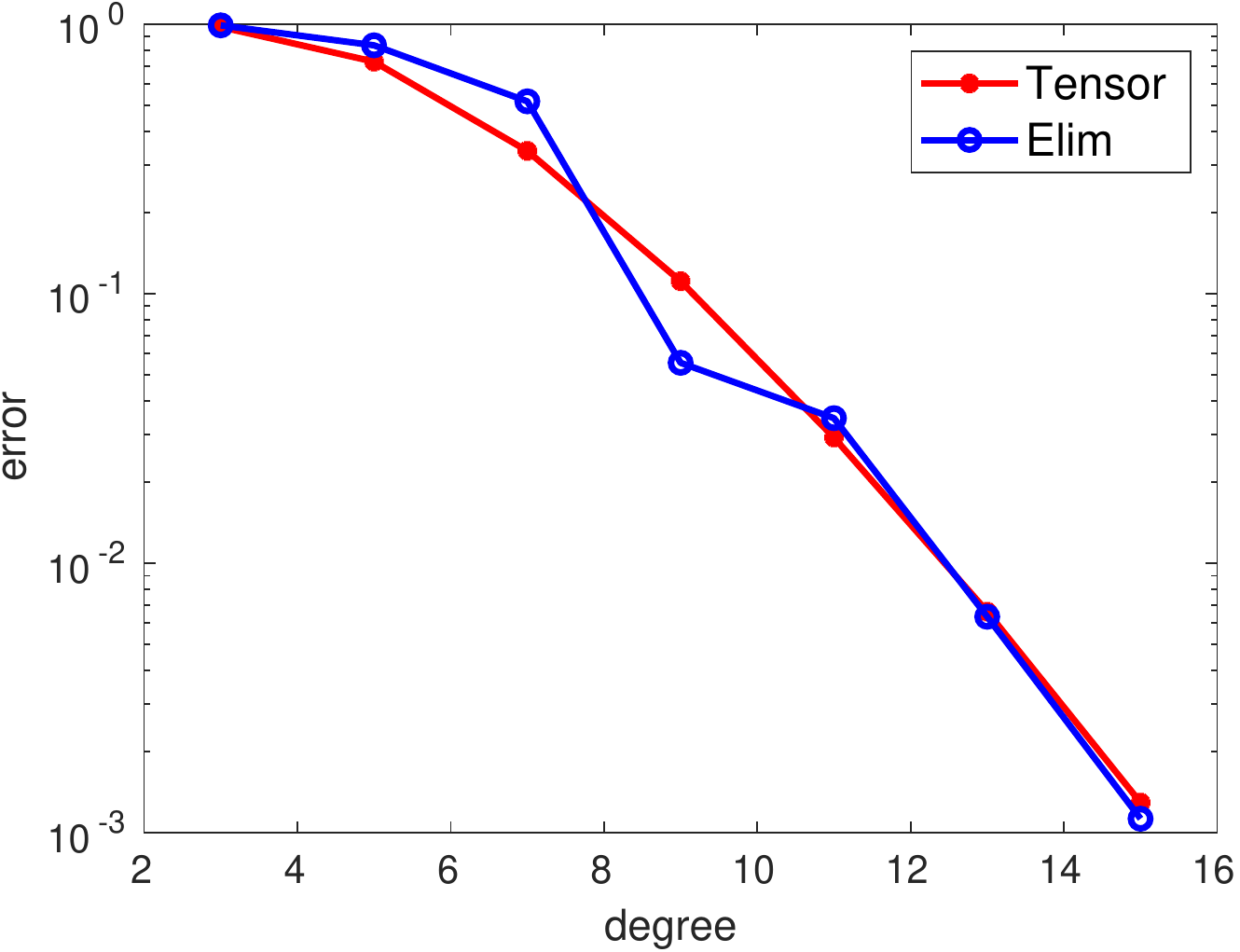}
\includegraphics[width=2.2in]{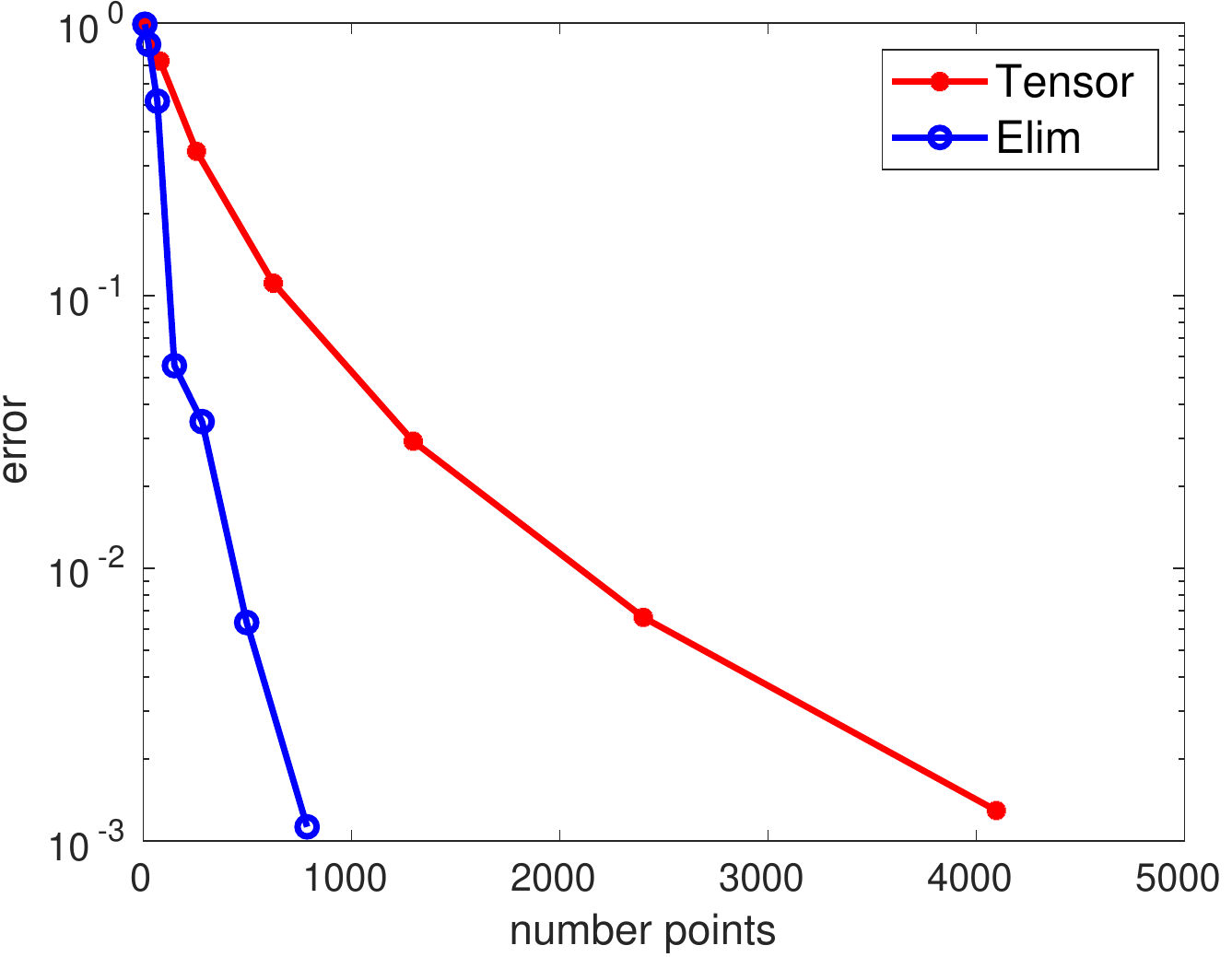}
\end{center}
\caption{Comparison of the tensor product rules 
  with the rules of table~\ref{tab:S4rules}. 
  $\Omega=T_4$. Left: error vs. degree. Right: error vs. number points.}
\label{fig:fourD}
\end{figure}

For higher dimensional simplices, the standard
approach is to use the Duffy transform and use tensor product Gauss-Jacobi rules.
The plots \ref{fig:threeD} and \ref{fig:fourD}
compare the quadrature error of these rules with the rules of tables~\ref{tab:S3rules}
and \ref{tab:S4rules}. 

One can see that for the simplices the accuracy of the tensor product
and the newly obtained rules of the same degree are very similar.
This is despite the fact that the tensor product rules integrate
polynomials in the space
$\{\xb^{\boldsymbol{\alpha}}: \alpha_i\leq p\}$ exactly, which is a
much larger space than $\Pol^d_p$.  Apparently, these additional
polynomials do not contribute to the accuracy in our examples.  The
right figure shows the same data, but here the error is plotted versus
the number of quadrature nodes. This illustrates that the rules
obtained with the elimination process
require much fewer function evaluations to achieve the
same accuracy. The gain of using these rules increases with the dimension.

The obvious approach for the domain $T_2\times T_2$ is to use the quadrature rule in
\eqref{quad:S2:Duffy} for both triangles. Alternatively, one can use
the rules of table~\ref{tab:S4rules} . The resulting
cubature errors are shown in figure~\ref{fig:fourD}. 

Here it is interesting that for the same degree the tensor product rule is
more accurate but the generalized Gauss rule is much more
accurate for the same number of cubature nodes.

\begin{figure}[hbt]
\begin{center}
\includegraphics[width=2.2in]{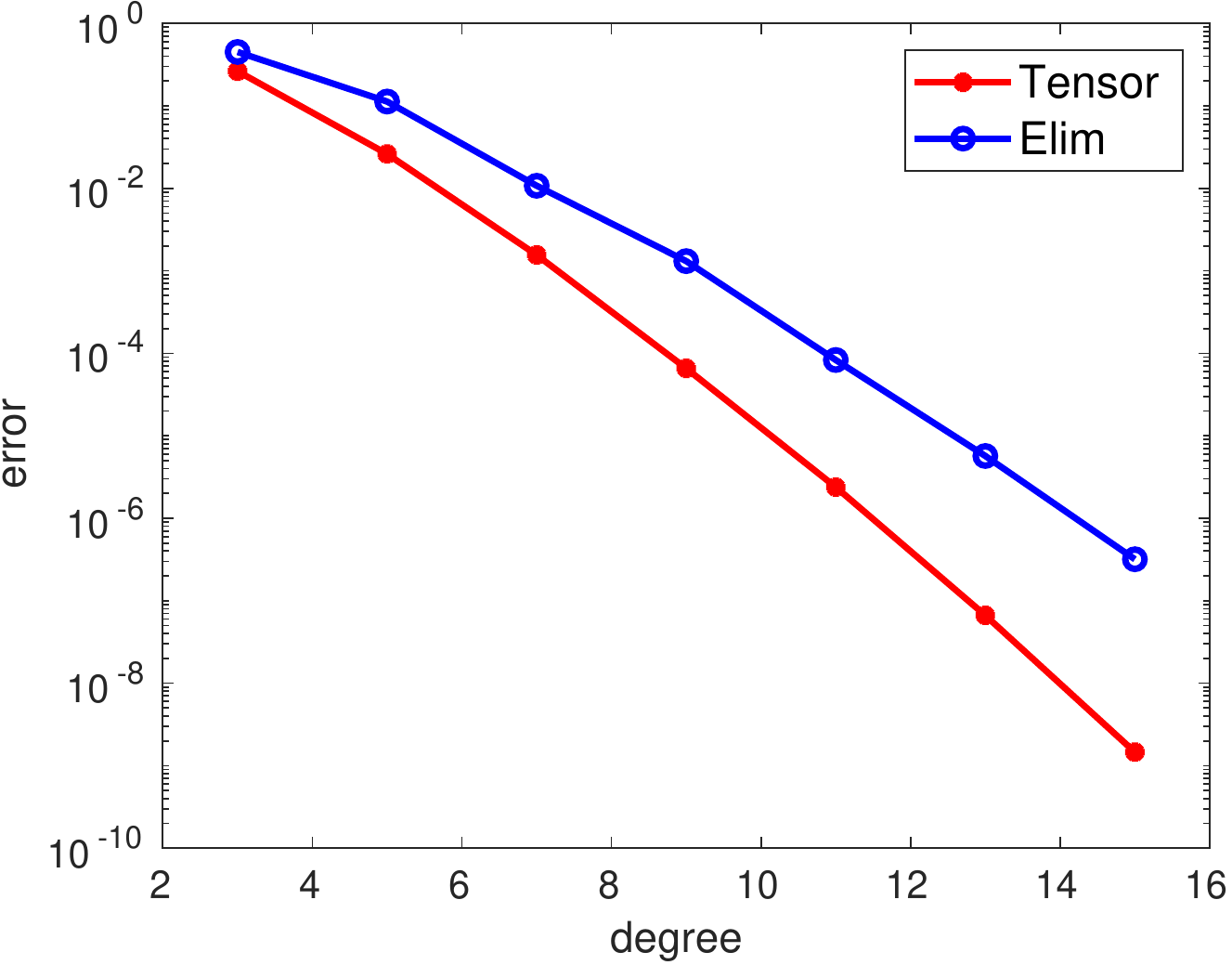}
\includegraphics[width=2.2in]{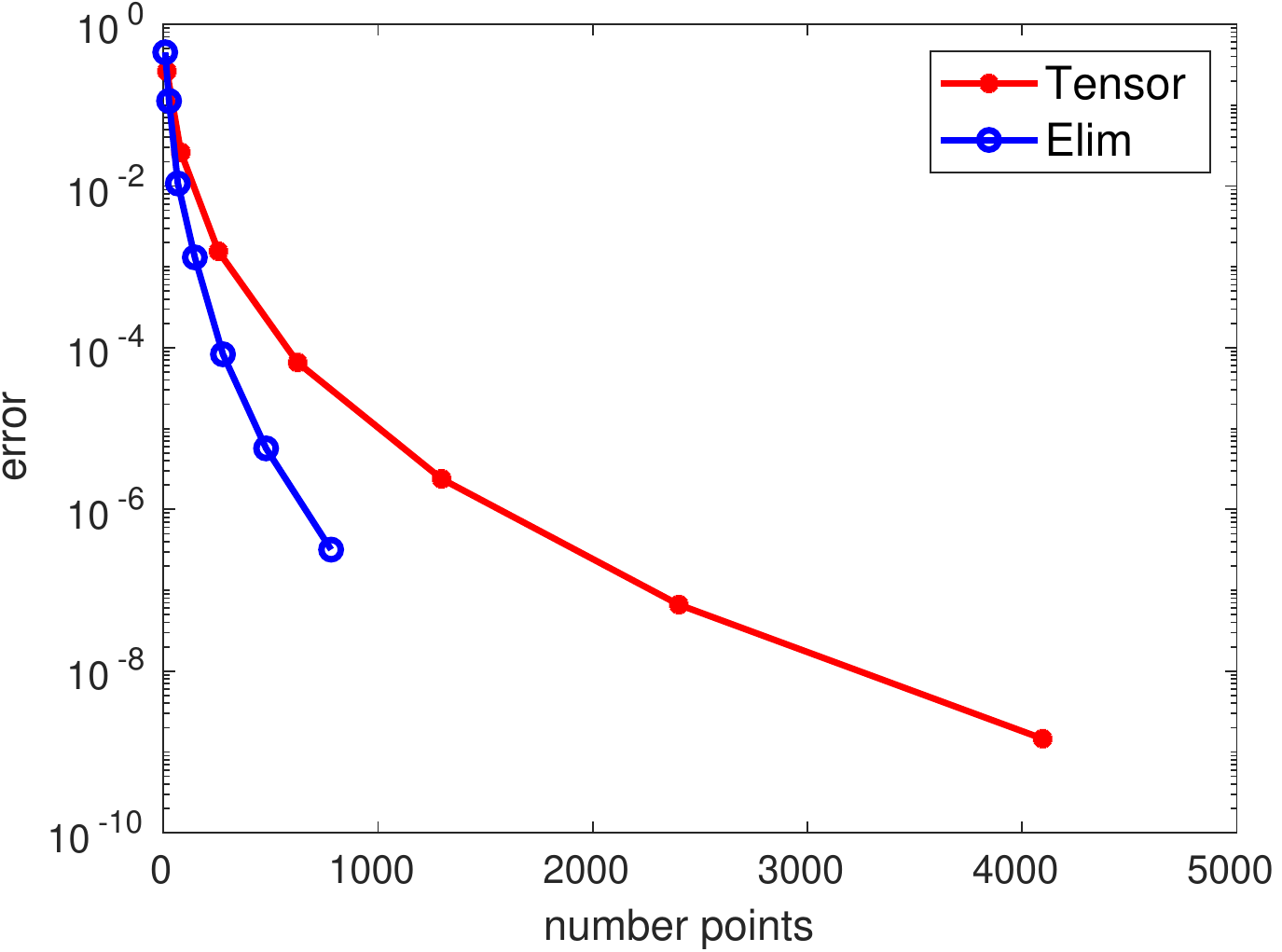}
\end{center}
\caption{Comparison of the tensor product rules 
  with the rules of table~\ref{tab:S2S2rules}.}
\label{fig:SimplexSimplex}
\end{figure}

\section{Conclusions}
We have developed a new node elimination scheme, implemented it for
several high-dimensional polytopes and obtained high-degree cubature
rules. In the node elimination the degree of precision of a cubature
rule is preserved while the number of nodes is reduced. For the
integrands considered here,
the node elimination does not reduce the accuracy in the case that the
domain is a simplex, while for $T_2\times T_2$ the
accuracy is somewhat reduced. However, in all cases the trade off between cost
in terms of evaluation points and accuracy is greatly improved with
the newly obtained rules. Why two cubature rules of the same degree
of precision can have different accuracies appears to be an open
research question.

It is well known that symmetries of the domain can be exploited to
reduce the number of unknowns of the moment equations. This is of
practical importance for domains in higher dimensions than what is
presented in this work. However, this also gets more difficult because
of the more complicated structure of the symmetry group of high
dimensional polytopes. This is another area where further development
of the method is needed.


\end{document}